%% file: haagerup_dual.tex
\pdfoutput=1
\documentclass[12pt]{amsart}

\usepackage[bibtex, julia]{ulthiel}
\usetikzlibrary{babel}

\AtEveryBibitem{%
    \clearfield{urlyear}%
}

\usepackage[utf8]{inputenc}
\usepackage[T1]{fontenc}
\usepackage{amssymb}
\usepackage{autobreak}
\usepackage[linesnumbered,ruled,vlined]{algorithm2e}
\usepackage{tabularray} 

\SetKwInput{KwInput}{Input}                
\SetKwInput{KwOutput}{Output}              

\usepackage{bbold}
\usepackage{enumitem}
\usepackage[]{tcolorbox}

\input{Tikzpictures}
\tcbuselibrary{listings, breakable}

\presetkeys%
    {todonotes}%
    {inline,backgroundcolor=yellow}{}

\usepackage{lstautogobble}
\usepackage{IEEEtrantools}

\usepackage{listofitems} 
\newcommand{\bbk}{\mathbb{k}}

\newcommand{\Nsd}[1]{
	\readlist*\Nvar{#1}%
	N_{{\Nvar[1]},{\Nvar[2]}}^{\Nvar[3]}
}

\newcommand{\Ft}[1]{%
	\readlist*\Fvar{#1}%
	[F_{{\Fvar[1]}, {\Fvar[2]}, {\Fvar[3]}}^{\Fvar[4]}]^{ ( {\Fvar[5]}, {\Fvar[6]}, {\Fvar[7]} ) }_{ ( {\Fvar[8]}, {\Fvar[9] }, {\Fvar[-1] } ) }
}

\newcommand{\Finvt}[1]{%
\readlist*\Fvar{#1}%
[(F_{{\Fvar[1]}, {\Fvar[2]}, {\Fvar[3]}}^{\Fvar[4]} )^{-1} ]^{ ( {\Fvar[5]}, {\Fvar[6]}, {\Fvar[7]} ) }_{ ( {\Fvar[8]}, {\Fvar[9] }, {\Fvar[-1] } ) }
}
\newcommand{\Fs}[1]{%
	\readlist*\Fvar{#1}%
	[F_{{\Fvar[1]}, {\Fvar[2]}, {\Fvar[3]}}^{\Fvar[4]}]^{{\Fvar[5]}}_{{\Fvar[6]}}
}

\newcommand{\Rs}[1]{%
	\readlist*\Fvar{#1}%
	R_{{\Fvar[1]}, {\Fvar[2]}}^{\Fvar[3]}
}

\newcommand{\Rt}[1]{%
	\readlist*\Rvar{#1}%
	[R_{{\Rvar[1]},{\Rvar[2]}}^{\Rvar[3]}]^{\Rvar[4]}_{\Rvar[5]}
}

\newcommand{\g}[1]{%
	\readlist*\gvar{#1}%
	g_{{\gvar[1]}, {\gvar[2]}}^{\gvar[3]}
}

\newcommand{\fpmod}{%
  \mathcal{M}_\Phi^{|L_{FP}^\Phi|}
}

\newcommand{\lfpp}{%
  L_{FP}^\Phi
}

\newcommand{\gaugespace}{ \left\{ g_{ij}^k \in \mathbb{C} \backslash \{0\} \right\} }

\newcommand{\calc}{ \mathcal C }

\newcommand{\HI}{\mathcal H _3}

\title[$F$-symbols and $R$-symbols for the center of the Haagerup subfactor]{$F$-symbols and $R$-symbols for the Drinfeld center of the Haagerup subfactor}

\author[F. Mäurer]{Fabian Mäurer$^*$}
\email{maeurer@mathematik.uni-kl.de}

\author[U. Thiel]{Ulrich Thiel$^*$}
\email{ulrich.thiel@math.rptu.de}

\address{$^*$Department of Mathematics, RPTU University Kaiserslautern-Landau, Postfach 3049, 67663 Kaiserslautern, Germany}

\author[G. Vercleyen]{Gert Vercleyen$^\dagger$}
\email{gvercley@purdue.edu}
\address{$^\dagger$Department of Mathematics, Purdue University, 150 N University St, West Lafayette, 47907, Indiana, USA}

\date{January 27, 2026}

\DeclareMathOperator{\coker}{coKer}

\DeclareUnicodeCharacter{2295}{\ensuremath{\oplus}}
\DeclareUnicodeCharacter{2297}{\ensuremath{\otimes}}
\DeclareUnicodeCharacter{1D7D9}{\ensuremath{\mathbb 1}}
\DeclareUnicodeCharacter{2124}{\ensuremath{\mathbb Z}}
\DeclareUnicodeCharacter{2083}{\ensuremath{{}_3}}
\DeclareUnicodeCharacter{1D49E}{\ensuremath{\mathscr C}}
\DeclareUnicodeCharacter{1D4B5}{\ensuremath{\mathscr Z}}
\DeclareUnicodeCharacter{1D49F}{\ensuremath{\mathscr D}}
\DeclareUnicodeCharacter{3C4}{\ensuremath{\tau}}
\DeclareUnicodeCharacter{3D5}{\ensuremath{\phi}}
\DeclareUnicodeCharacter{221A}{\ensuremath{\sqrt{}}}
\DeclareUnicodeCharacter{3C7}{\ensuremath{\chi}}
\DeclareUnicodeCharacter{22C5}{\ensuremath{\cdot}}
\DeclareUnicodeCharacter{D7}{\ensuremath{\times}}
\DeclareUnicodeCharacter{3B1}{\ensuremath{\alpha}}
\DeclareUnicodeCharacter{3B3}{\ensuremath{\gamma}}
\DeclareUnicodeCharacter{3C1}{\ensuremath{\rho}}
\DeclareUnicodeCharacter{2217}{\ensuremath{\ast}}
\DeclareUnicodeCharacter{22EE}{\ensuremath{\vdots}}

\begin{document}

\begin{abstract}
We have recently devised and implemented an algorithm for computing the Drinfeld center of a fusion category in our software package \textsc{TensorCategories.jl}. By simplifying the field of definition of the $F$-symbols, $R$-symbols, and pivotal coefficients of the Haagerup $\mathcal{H}_3$ category, we were now able to compute explicit $F$-symbols, $R$-symbols, and pivotal coefficients for its Drinfeld center $\mathcal{Z}(\HI)$. In this paper, we discuss how to access this data. We also present our new general algorithm for finding a gauge of a multiplicity-free fusion category in which the symbols belong to a minimal field.
\end{abstract}

\maketitle

\newcommand{\tageq}{\refstepcounter{equation}\tag{\theequation}}

\newcommand{\mat}[1]{\begin{bmatrix}
    #1
\end{bmatrix}}

\section{Introduction}
Fusion categories appear in a wide variety of branches in mathematics and physics such as
topological quantum computation \cite{freedmanTopologicalQuantumComputation2002,kitaevAnyonsExactlySolved2006,wang2010topological,rowellInvitationMathematicsTopological2016,rowellMathematicsTopologicalQuantum2018},
knot theory \cite{kasselQuantumGroups1995,reshetikhinRibbonGraphsTheir1990,reshetikhinInvariants3manifoldsLink1991,turaevQuantumInvariantsKnots2010},
representation theory of weak Hopf algebras and quantum groups \cite{kasselQuantumGroups1995,bakalovLecturesTensorCategories2001},
low-dimensional topology and topological field theory \cite{turaevMODULARCATEGORIES3MANIFOLD1992,craneStateSumInvariants4Manifolds1997,bakalovLecturesTensorCategories2001,cuiStateSumInvariants2017,turaevMonoidalCategoriesTopological2017},
subfactor theory and planar algebras \cite{jonesClassificationSubfactorsIndex2013,grossmanExtendedHaagerupFusion2018,calegariCyclotomicIntegersFusion2011,evansNeargroupFusionCategories2012},
vertex operator algebras and conformal field theory \cite{mooreClassicalQuantumConformal1989,creutzigFusionCategoriesAffine2019,evansNeargroupFusionCategories2012}.

For any practical approach for working with fusion categories, the skeletal data, i.e., the coefficients of the linear maps encoding the associator, pivotal structure, and, if applicable, braiding, should be known. These coefficients, which we will call $F$-symbols, $P$-symbols, and $R$-symbols respectively, could, at least theoretically, be obtained from a fusion ring of the corresponding fusion categories by solving systems of polynomial equations and inequalities. This process is called categorifying the fusion ring. In \cite{osborne2019f}, the skeletal data of the $\mathcal{H}_3$ Haagerup category was obtained via categorification of its ring, and later in \cite{vercleyen2024lowrankmultiplicityfreefusioncategories}, this approach has been used to find all multiplicity-free fusion categories with braided and pivotal structures up to rank seven.
These are accessible on the AnyonWiki \cite{anyonwiki}, as part of the Anyonica software package \cite{anyonica}, and, recently, also as part of the \textsc{TensorCategories.jl} package \cite{Maeurer_TensorCategories_jl}.

Unfortunately, categorifying fusion rings by solving systems of polynomial equations has only proven fruitful for multiplicity-free fusion rings of low rank. This is because the number of polynomial equations rises dramatically with multiplicity and the Frobenius-Perron dimension of a fusion ring.
An alternative way to obtain examples of categories is to apply certain constructions to smaller categories to obtain larger ones. The specific construction this paper concerns is the Drinfeld center. \\

Recently, we have devised and implemented an algorithm to compute the $F$-symbols, $R$-symbols, and $P$-symbols (called symbols from now on) of the Drinfeld center $\mathcal Z(\mathcal C)$ of a pivotal fusion category $\mathcal C$ as part of our software package \textsc{TensorCategories.jl} \cite{center-computation,Maeurer_TensorCategories_jl}.
This is an open source software package based on the programming language Julia and using the OSCAR \cite{OSCAR,OSCAR-book} computer algebra system for its symbolic computations.
By using the data from \cite{vercleyen2024lowrankmultiplicityfreefusioncategories}, the software was already used to compute the symbols of the centers of all multiplicity-free fusion categories up to rank $5$.

The speed with which the algorithm computes symbols of $\mathcal Z (\mathcal C)$ depends, among other factors, on the field the symbols of $\mathcal C$ belong to (a subfield of $\mathbb{C}$). Ideally, the symbols are expressed in a cyclotomic field. For almost all multiplicity-free fusion categories up to rank~$7$, this is indeed possible. 

However, there are $4$ exceptions \cite{morrison2012non}: two unitary categories
$\text{HI}(\mathbb{Z}_3)_{1,0,1} \cong \HI$ and
$\text{HI}(\mathbb{Z}_3)_{2,0,1} \cong \mathcal{H}_2$, and two non-unitary categories
$\text{HI}(\mathbb{Z}_3)_{3,0,1}$ and
$\text{HI}(\mathbb{Z}_3)_{4,0,1}$. 
These categories are, up to Morita equivalence, the possible categorifications of the fusion ring in Table \ref{fig:H3_fusion}, and we use here the notation from the AnyonWiki to denote the specific categorifications. The categories $\HI$ and $\mathcal H_2$ are related via a Galois conjugation.
\begin{table}[htbp]
    \caption{The fusion rules for $\mathcal H_2$ and $\mathcal H_3$.}
    \label{fig:H3_fusion}
    \begin{tabular}{c||c|c|c|c|c|c}
        & $\mathbb 1$ & $\alpha$ & $\alpha_\ast$ & $\rho$ & $\alpha\rho$ & $\alpha_\ast\rho$ \\ \hline \hline
        $\mathbb 1$ & $\mathbb 1$ & $\alpha$ & $\alpha_\ast$ & $\rho$ & $\alpha\rho$ & $\alpha_\ast\rho$ \\ \hline
        $\alpha$ & $\alpha$ & $\alpha_\ast$ & $\mathbb 1$ & $\alpha \rho$ & $\alpha_\ast\rho$ & $\rho$ \\ \hline
        $\alpha_\ast$ & $\alpha_\ast$ & $\mathbb 1$ & $\alpha$ & $\alpha_\ast\rho$ & $\rho$ & $\alpha\rho$ \\ \hline
        $\rho$ & $\rho$ & $\alpha_\ast\rho$ & $\alpha\rho$ & $\mathbb 1 + \rho + \alpha\rho + \alpha_\ast\rho$ & $\alpha_\ast + \rho + \alpha\rho + \alpha_\ast\rho$ & $\alpha + \rho + \alpha\rho + \alpha_\ast\rho$ \\ \hline
        $\alpha\rho$ & $\alpha\rho$ & $\rho$ & $\alpha_\ast\rho$ & $\alpha + \rho + \alpha\rho + \alpha_\ast\rho$ & $\mathbb 1 + \rho + \alpha\rho + \alpha_\ast\rho$ & $\alpha_\ast + \rho + \alpha\rho + \alpha_\ast\rho$ \\ \hline
        $\alpha_\ast\rho$ & $\alpha_\ast\rho$ & $\alpha\rho$ & $\rho$ & $\alpha_\ast + \rho + \alpha\rho + \alpha_\ast\rho$ & $\alpha + \rho + \alpha\rho + \alpha_\ast\rho$ & $\mathbb 1 + \rho + \alpha\rho + \alpha_\ast\rho$
    \end{tabular}
\end{table}
These four categories stem from the Haagerup subfactor \cite{haagerup1994principal} and are often called \emph{exotic} in the sense that they do not seem to be related to any known construction from finite groups and/or quantum groups. \\

Categories coming from the Haagerup subfactor are not just a mathematical curiosity. In the theory of rational conformal field theories (RCFTs), modular fusion categories (MFCs) appear as representation categories of the vertex operator algebra of an RCFT \cite{huang2005vertex,huang2008rigidity,EGNO}. An open question is whether any modular fusion category leads to an associated RCFT. This seems to be the case for many (MFCs) \cite{huang2014logarithmic}, but there has been no conclusive answer for Drinfeld centers of categories coming from the Haagerup subfactor yet, despite recent numerical and analytic efforts \cite{huang2022numerical,vanhove2022critical}. The exoticness of the $\mathcal{ H}_3$ category has gained considerable interest from condensed matter theory where lattice  $\mathcal{H}_3$ symmetry are being investigated\cite{bottini2025construction,jiaH3,corcoranH3,liuH3}. The central object of interest in these attempts has been $\mathcal Z ( \mathcal H_3)$. We hope that our results will help to finally answer this open question. \\

The Drinfeld center $\mathcal{Z}(\mathcal H_3)$ is an enormous object: it has over 1.6 million $F$-symbols and 1400 $R$-symbols. In our first attempt at computing the center, we used the $F$-symbols for $\mathcal H_3$ from \cite{osborne2019f}, \cite{huang2020f}, and \cite{vercleyen2024lowrankmultiplicityfreefusioncategories}, which are contained in a number field $K$ of degree 8. While we could indeed compute the center over $K$, it did not split there, and we were not able to pass to a splitting field, i.e., we did not manage to compute the center over $\mathbb{C}$. We discussed the technicalities of splitting in \cite{center-computation}. 

In this paper, we present a general algorithm that expresses symbols of a multiplicity-free fusion category $\mathcal{C}$ in a gauge that requires a \emph{minimal field} $M(\mathcal{C})$, see section \ref{sec:fields}. Using this algorithm, we managed to express the symbols of $\HI$ over a number field of degree $4$ as compared to $8$ previously. Specifically:

$$
\begin{aligned}
  & M(\mathcal{H}_2) = \mathbb{Q}(\alpha) \;, \quad \alpha^4 - \alpha^3 - \alpha^2 - \alpha + 1 = 0 \;, \\
  & M(\mathcal{H}_3) = \mathbb{Q}(\beta) \;, \quad \beta^4 + \beta^3 - \beta^2 + \beta + 1 = 0 \;.
\end{aligned}
$$

Using these simplified symbols, we were then indeed able to compute all symbols of $\mathcal{Z}(\mathcal H_3)$ over a splitting field. Our data is freely available (both numerically and symbolically) as part of the \textsc{TensorCategories.jl} package and as part of an online repository \cite{ZH3_Data_Repo}. The computation of the simple objects took about 100 seconds, and the computation of the $F$-symbols took about 3 days.

In \ref{sec:code} we show how this data can be accessed. In section \ref{sec:fields}, we provide our new algorithm to express the symbols of a multiplicity-free pivotal fusion category in a gauge for which the smallest field over which they are defined is a minimal field for the fusion category. We begin in section \ref{sec:Preliminaries} with a review of the definitions and conventions on fusion categories, their representation via fusion systems, and their center.


\subsection*{Acknowledgements}
This work is a contribution to the SFB-TRR 195 ``Symbolic Tools in
Mathematics and their Application'' of the German Research Foundation (DFG), Project-ID 28623755.

\newpage

\tableofcontents

\section{Preliminaries}\label{sec:Preliminaries}

We will recall the central definitions used and refer to \cite{EGNO,center-computation} for anything further.

\subsection{Fusion Categories and Fusion Systems}\label{sec:fuscat}

Let $\mathbb k$ be a field. We do not require $\mathbb k$ to be algebraically closed.
\begin{definition}
    A \emph{tensor category} is a locally finite $\mathbb k$-linear abelian rigid monoidal category such that the unit is simple and the monoidal product is bilinear on morphisms. A \emph{fusion category} of \emph{rank} $r$ is a semisimple tensor category with finitely many non-isomorphic simple objects $X_1,...,X_r$ such that $\End(X_i) \cong \mathbb k$ for all $i = 1,...,r$.
\end{definition}

For all details on tensor categories over arbitrary fields, we refer to \cite{center-computation} or \cite{bruguieres2013center}.
Let $\mathcal C$ be a fusion category over $\bbk$ with simple objects $X_1,...,X_r$. We call the spaces $H_{i,j}^k := \Hom(X_i \otimes X_j, X_k)$ the \emph{multiplicity spaces}. The naming makes sense since for any fusion category $\dim_{\bbk}H_{i,j}^k$ is the multiplicity of $X_k$ in $X_i\otimes X_j$. The maximum $\max\{\dim H_{i,j}^k\}$ is called the \emph{multiplicity} of the category. We fix a basis for all multiplicity spaces $H_{i,j}^k$. The freedom in this choice is known as \emph{gauge freedom} and choosing the basis is also called \emph{fixing the gauge}.

The choice of bases for the multiplicity spaces gives rise to canonical bases of the spaces $$\Hom((X_i \otimes X_j) \otimes X_k, X_l),~~\text{and}~~\Hom(X_i \otimes (X_j \otimes X_k), X_l)$$ by factoring through all the simple $X_m$ and $X_n$ like this:

\begin{align}
    \tzCanonicalBasisLeft{X_i}{X_j}{X_k}{X_l}{X_m}{\alpha}{\beta}~~\text{and}~~\tzCanonicalBasisRight{X_i}{X_j}{X_k}{X_l}{X_n}{\gamma}{\delta}
\end{align}
where $\alpha \in \Hom(X_i\otimes X_j, X_m)$, $\beta \in \Hom(X_m \otimes X_k, X_l)$, $\gamma \in \Hom(X_j \otimes X_k, X_n)$ and $\Hom(X_i \otimes X_n, X_l)$.

The associator morphisms $a_{X,Y,Z} \colon (X \otimes Y) \otimes Z \to X \otimes (Y \otimes Z)$ are then encoded by the matrices representing the induced map

\begin{align}
    \begin{aligned} \label{F_matrix_iso}
      \Hom((X_i \otimes X_j) \otimes X_k, X_l) &\simeq \Hom(X_i \otimes (X_j \otimes X_k), X_l) \\
      f &\mapsto f \circ a_{X_i,X_j,X_k}^{-1}.
      \end{aligned}
\end{align}

In particular, the image of this map induced by the associator acting on basis vectors can be expressed as
\begin{align}
    \tzCanonicalBasisLeft{X_i}{X_j}{X_k}{X_l}{X_m}{\alpha}{\beta} \circ a_{X_i,X_j,X_k}^{-1} = \sum\limits_{n = 1}^r\sum\limits_{\substack{\gamma \in H_{j,k}^n\\\delta \in H_{i,n}^l}} \left[F_{i,j,k}^l\right]_{(n,\gamma,\delta)}^{(m,\beta,\alpha)}\tzCanonicalBasisRight{X_i}{X_j}{X_k}{X_l}{X_n}{\gamma}{\delta}
\end{align}

where $\Ft{i,j,k,l,n,\beta,\alpha,m,\gamma,\delta} \in \bbk$. The coefficients, $\Ft{i,j,k,l,n,\beta,\alpha,m,\gamma,\delta}$, are known as \emph{$F$-symbols} of $\mathcal C$. Note that the convention for the order of the labels is not consistent in the literature. In all our computations, we will use this precise ordering, and also the exported data uses the same.

If $\mathcal C$ is braided with braided structure $\sigma_{X,Y} \colon X \otimes Y \simeq Y \otimes X$ the same procedure leads to a representation of the braiding in terms of the chosen bases, i.e., the image of the braiding under the isomorphism
\begin{align}
    \Hom(X_i \otimes X_j, X_k) &\simeq \Hom(X_j \otimes X_i, X_k) \\
    f &\mapsto f \circ \sigma_{X_i, X_j}^{-1}
\end{align}
can be expressed in terms of matrices as
\begin{align}
    \tzBasisVector{X_i}{X_j}{X_k}{\alpha} \circ \sigma_{X_i,X_j}^{-1} = \sum\limits_{\beta \in H_{i,j}^k}^r \left[R_{i,j}^k\right]^{\alpha}_\beta \tzBasisVector{X_i}{X_j}{X_k}{\beta}
\end{align}
where the values $\left[R_{i,j}^k\right]^{\alpha}_\beta \in \bbk$ are called \emph{$R$-symbols}.

Likewise, if $\mathcal C$ is pivotal with pivotal structure $\psi_{X}: X \rightarrow X^{**}$, then a gauge choice leads to a representation of $\psi$ as a list $\left\{P_i \in \bbk | i = 1, \ldots, r \right\}$ called the \emph{pivotal coefficients} or $P$-symbols of $\calc$.

In practice, a fusion category $\mathcal{C}$ over a field $\bbk$ is often represented as a fusion system, which is defined as follows
\begin{definition} Let $\mathcal{C}$ be a fusion category over $\bbk$. Via a gauge choice, $\mathcal{C}$ can be represented as a \emph{fusion system}, which is a tuple $\Phi_\calc =  ( \mathbf{L},*,\mathbf{N}, \mathbf{F}) $  where
  \begin{enumerate}
    \item  $\mathbf{L} = \{ 1, \ldots, r\}$ is a tuple of representatives (or labels) for the isomorphism classes of simple objects of $\mathcal{C}$. Here, $r$ is the rank of the fusion category.
    \item  $*:\mathbf{L} \rightarrow \mathbf{L}$ is a map for which $1^* = 1$ and $(i^* )^* = i, \forall i \in \mathbf{L}$, which comes from the rigidity of the $\mathcal{C}$
    \item  $\mathbf{N} = \left\{ \left. N_{i,j}^k \right| i,j,k = 1, \ldots,r \right\}  $, with $N_{i,j}^k :=\text{dim}(H^k_{i,j})$,  is the set of structure constants of the Grothendieck ring of $\mathcal{C}$, and
    \item  $\mathbf{F} =  \left\{ \Ft{i,j,k,l,n,\delta,\gamma,m,\beta,\alpha} \in \bbk \right\} $ is the list of $F$-symbols representing the associator. Their values depend on the specific gauge choice.
  \end{enumerate}
  If $\calc$ is multiplicity-free, i.e. $\max\dim\{H_{i,j}^k\} = 1$ then we call $\Phi_\calc$ a \emph{multiplicity-free fusion system}.
\end{definition}
By adding extra structure(s), we obtain the notions of pivotal, spherical, braided, ribbon, and modular fusion systems.
\begin{definition}
Let $\mathcal{C}$ be a fusion category over $\bbk$.
\begin{itemize}
  \item If $\mathcal{C}$ is a pivotal fusion category over $\bbk$, then, via a choice of gauge, it can be represented as a \emph{pivotal fusion system} $( \mathbf{L} , *, \mathbf{N} , \mathbf{F} , \mathbf{P}  )$ which is a fusion system together with a list $\left\{P_i \in \bbk | i = 1, \ldots, r \right\}$ that represents the pivotal structure of $\mathcal C$
  \item The \emph{left quantum dimensions} $\{d_a^L \in \mathbb{C}| a \in L \}$ of a pivotal fusion system are defined as
  \begin{eqnarray}
    d_a^L := \frac{P_a}{\Ft{a^*,a,a^*, a^*,1,1,1,1,1,1}}.
  \end{eqnarray}
  The pivotal structure $\mathbf{P}$ is called \emph{spherical} if $d_a^L = d^L_{a^*}$ for all $a \in \mathbf{L}$. In this case, the pivotal fusion system is called a \emph{spherical fusion system}, the left quantum dimensions are called the \emph{quantum dimensions} of the system and are denoted by $d_a$ rather than $d_a^L$.
  \item If $\mathcal{C}$ is a braided fusion category over $\bbk $, then, via a choice of gauge, it can be represented as a \emph{braided fusion system} $( \mathbf{L} , *, \mathbf{N} , \mathbf{F} , \mathbf{R}  )$ which is a fusion system together with a finite list of $R$-symbols
  $\left\{\Rt{a,b,c,\alpha,\beta} \in \bbk \right\} $ that represent the braided structure.
  \item A \emph{ribbon fusion system} $( \mathbf{L}, *, \mathbf{N}, \mathbf{F}, \mathbf{P}, \mathbf{R} )$ is a spherical braided fusion system.
  \item A \emph{modular fusion system} $( \mathbf{L}, *, \mathbf{N}, \mathbf{F}, \mathbf{P}, \mathbf{R} )$ is a ribbon fusion system for which
    the $S$-matrix matrix $ S \in \text{Mat}_{r \times r}(\mathbb{k})$, whose entries are given by
    \begin{eqnarray}
      [S]^a_b =
        \sum_{c = 1}^r
        \sum_{i = 1}^{\Nsd{a,b^*,c}}
        \sum_{j = 1}^{\Nsd{c,b,a}}
        \sum_{k = 1}^{\Nsd{b^*,a,c}}
        \sum_{l = 1}^{\Nsd{a,b^*,c}}
        d_a d_b
    \Finvt{a,b^*,b,a,1,1,1,c,j,i} \Rt{b^*,a,c,i,k} \Rt{a,b^*,c,k,l} \Ft{a,b^*,b,a,c,j,l,1,1,1},
    \end{eqnarray}
    is invertible.
\end{itemize}
\end{definition}

In \cite{davidovichArithmeticModularCategories2013}, it is shown that each of the systems above defines a fusion category, with the corresponding adjectives, that is unique up to equivalence.

However, because of the gauge freedom and fusion ring automorphisms, there are an infinite number of fusion systems that represent the same category. For multiplicity-free fusion categories $\dim H^k_{i,j} \in \{0,1\}$, so changing the bases of $H^k_{i,j}$ (i.e., applying a gauge transform) comes down to multiplying the previously chosen basis vector of each $H^k_{i,j}$ by a non-zero factor $\g{i,j,k}$ which belongs to a (possibly trivial) field extension of $\bbk$. This induces the following map $\varphi$ on the $F$-, $R$-, and $P$-symbols of the corresponding multiplicity-free fusion system.
\begin{align*}\label{eq:gt}
  \varphi:& \Fs{i,j,k,l,m,n}  \mapsto \frac{ \g{i,j,m} \g{m,k,l} }{ \g{i,n,l} \g{j,k,n} } \Fs{i,j,k,l,m,n},\\
  \varphi:& \Rs{i,j,k}        \mapsto \frac{ \g{i,j,k} }{ \g{j,i,k} } \Rs{i,j,k},\\
  \varphi:& P_i               \mapsto \frac{\g{i^*,i,1}\g{1,i^*,i^*}}{\g{i^*,1,i^*}\g{i,i^*,1}} P_i
\end{align*}
where we denoted $\Ft{i,j,k,l,m,1,1,n,1,1}$ as $\Fs{i,j,k,l,m,n}$ and $\Rt{i,j,k,1,1}$ as $\Rs{i,j,k}$.
A fusion ring automorphism $\sigma$ is a map $\mathbf{L}\to \mathbf L $ such that $N_{\sigma(i),\sigma(j)}^{\sigma(k)} = N_{i,j}^k$ for all $i,j,k \in \mathbf{L}$. For a multiplicity-free fusion category, an automorphism $\sigma$ induces a map on the $F$-, $R$-and $P$-symbols as follows
\begin{align*}
  \sigma:& \Fs{i,j,k,l,m,n} \mapsto \Fs{\sigma(i),\sigma(j),\sigma(k),\sigma(l),\sigma(m),\sigma(n)}\\
  \sigma:& \Rs{i,j,k} \mapsto \Rs{\sigma(i),\sigma(j),\sigma(k)}\\
  \sigma:& P_i \mapsto P_{\sigma(i)}
\end{align*}

The $\mathcal{H}_3$ category is a pivotal fusion category. Two multiplicity-free pivotal fusion systems, with the same structure constants $\mathbf{N}$, define the same fusion category, up to isomorphism, if and only if their $F$-and $P$-symbols are related by a combination of a gauge-transform and a fusion ring automorphism. In Section \ref{sec:fields} we will see that we can use this freedom to find a gauge for a fusion system of $\mathcal{H}_3$ for which the $F$-and $P$-symbols belong to a \emph{minimal field} which is necessary to reduce the computational complexity of finding the center $\mathcal{Z}(\mathcal{H}_3)$.

\subsection{Complex extensions of fusion categories}\label{sec:embedding}
The \textsc{TensorCategories.jl} package allows one to work with fusion categories over a variety of fields, ranging from abstract fields to embedded subfields of $\mathbb C$ to the algebraic closure of $\mathbb Q$ (and even over finite fields, but this lies outside the scope of the current paper). Working with an abstract field, rather than an embedded one, yields more flexibility and speed for computations. For applications to physics, however, complex fusion categories, i.e., over an embedded subfield of $\mathbb C$, are required. This is because complex conjugation of the $F$-, $R$-, and $P$-symbols plays a fundamental role in determining the unitarity of operators and the probabilistic outcomes of processes involving anyons. In this section, we revise the theory of complex extensions of fusion categories and show that a fusion category over an abstract field can lead to multiple inequivalent fusion categories over embedded fields.

Let $\mathcal C$ and $\mathcal D$ be $\mathbb k$-linear categories. The $\mathbb k$-linear Deligne product $\mathcal C \boxtimes_{\mathbb k} \mathcal D$ is the cocompletion of the category with objects $X \boxtimes Y$ and $\Hom(X \boxtimes Y, X' \boxtimes Y') = \Hom_{\mathcal C}(X,X') \otimes_{\mathbb k} \Hom_{\mathcal D}(Y,Y')$. If $\mathcal C$ and $\mathcal D$ are tensor categories, then $\mathcal C \boxtimes_{\mathbb k} \mathcal D$ is a tensor category as well. More details can be found in \cite{EGNO}.

\begin{definition}
  Let $\mathcal C$ be a fusion category over a field $\mathbb k$. Let $\mathbb k \hookrightarrow \mathbb K$ be a field extension and $\mathrm{Vec}_{\mathbb K}$ the category of finite-dimensional vector spaces over $\mathbb K$. Then the category $\mathcal C \boxtimes_{\mathbb k} \mathbb K \coloneq \mathcal C \boxtimes_{\mathbb k} \mathrm{Vec}_{\mathbb K}$ is called the \emph{extension of scalars} of $\mathcal C$ by $\mathbb K$.
\end{definition}

The category $\mathcal C \boxtimes_{\mathbb k} \mathbb K$ is a tensor category over $\mathbb K$. If $\mathbb k \hookrightarrow \mathbb K$ is a separable extension, then $\mathcal C \boxtimes_{\mathbb k} \mathbb K$ is also semisimple \cite{non-split-fusion,non-split-fusion-JA}. Moreover, if $\mathcal C$ is pivotal (spherical), then so is $\mathcal C \boxtimes_{\mathbb k} \mathbb K$. This implies that if $\mathcal C$ is fusion then so is $\mathcal C \boxtimes_{\mathbb k} \mathbb K$ \cite{EGNO,non-split-fusion,non-split-fusion-JA}.

Let $(\mathbf{L},\ast,\mathbf N, \mathbf F)$ be a fusion system corresponding to a fusion category $\mathcal C$ over $\mathbb C$. Then the values of the $F$-symbols $\mathbf F$ lie in a number field of $\mathbb Q \subset \mathbb k \subset \mathbb C$ of finite degree over $\mathbb Q$, since they form a solution to a system of rational polynomial equations. Therefore. we can consider the fusion system $(\mathbf{L},\ast,\mathbf N, \mathbf F')$ where $F'$ is the restriction of $\mathbf F$ to $\mathbb k$. Then the corresponding fusion category $\mathcal C'$ is defined over $\mathbb k$ and $\mathcal C \cong \mathcal C' \boxtimes_{\mathbb k} \mathbb C$.

Given a fusion system (or category) defined over an abstract field $\mathbb k$, a complex fusion category can be obtained by choosing an embedding $\mathbb k \hookrightarrow \mathbb C$ and considering the extension of scalars $\mathcal C \boxtimes_{\mathbb k} \mathbb C$.

\begin{definition}
  The category $\mathcal C \boxtimes_{\mathbb k} \mathbb C$ is called a \emph{complex extension} of $\mathcal C$.
\end{definition}

Note that different complex extensions can lead to non-equivalent complex fusion categories.
\begin{example}
The two sets of $F$-symbols that solve the pentagon equations for the Fibonacci fusion ring are algebraically indistinguishable. However, the two Galois conjugated embeddings of $\mathbb Q(\phi)$ where $\phi^2 - \phi - 1 = 0$ lead to non-equivalent complex fusion categories. For $\phi = \frac{1+\sqrt{5}}{2}$ this is the Fibonacci category and  for $\phi = \frac{1-\sqrt{5}}{2}$ this is the Yang-Lee category. These two categories are not equivalent, since the first is unitary while the second is not.
\end{example}

\subsection{The Drinfeld center of a monoidal category}

Let $\mathcal C$ be a monoidal category.

\begin{definition}
  Let $X \in \mathcal C$. A \emph{half-braiding} for $X$ is a natural isomorphism
  \begin{align}
    \gamma = \{\gamma(Y) \colon X \otimes Y \to Y \otimes X \mid Y \in \mathcal C\}
  \end{align}
  such that
  \begin{equation}
    \begin{tikzcd}
      & (Y \otimes X) \otimes Z \ar[r, "a_{Y,X,Z}"] & Y \otimes (X \otimes Z) \ar[rd, "\id_Y \otimes \gamma(Z)"] &  \\
      (X \otimes Y) \otimes Z \ar[ur, "\gamma(Y) \otimes \id_Z"] \ar[dr, "a_{X,Y,Z}"] & & & Y \otimes (Z \otimes X) \\
      & X \otimes (Y \otimes Z) \ar[r, "\gamma(Y \otimes Z)"] & (Y \otimes Z) \otimes X \ar[ru, "a_{Y,Z,X}"] &
    \end{tikzcd}
  \end{equation}
  commutes for all $Y,Z \in \mathcal C$ and $\gamma(\mathbb 1) = \id_X$.
\end{definition}

\begin{definition}
  The \emph{Drinfeld center} (or just \emph{center}) $\mathcal Z(\mathcal C)$ of $\mathcal C$ consists of tuples $(X,\gamma)$ where $X \in \mathcal C$ and $\gamma$ is a half-braiding. Morphisms between tuples $(X, \gamma_X)$ and $(Y, \gamma_Y)$ are given by morphisms $f \colon X \to Y$ such that
  \begin{equation}
    \begin{tikzcd}
      X \otimes Z \ar[r, "f \otimes \id_Z"] \ar[d, "\gamma_X(Z)"] & Y \otimes Z \ar[d, "\gamma_Y(Z)"] \\
      Z \otimes X \ar[r, "\id_Z \otimes f"] & Z \otimes Y
    \end{tikzcd}
  \end{equation}
  commutes for all $Z \in \mathcal C$.
\end{definition}


%
%

\section{The data of $\mathcal Z(\mathcal H_3)$}\label{sec:code}

In this section, we explain how one can obtain all the data we computed. We assume that the reader has access to a computer with the latest version of Julia and the Oscar and \textsc{TensorCategories.jl} packages installed as described in \cite{center-computation,Maeurer_TensorCategories_jl}. The first step after starting Julia is to load the package as follows
\begin{tcolorbox}[title = \juliasymbol, breakable]
    \begin{juliaprompt}
      julia> using TensorCategories, Oscar
    \end{juliaprompt}
\end{tcolorbox}

The center $\mathcal{Z}(\HI)$ category can be accessed as follows.
\begin{tcolorbox}[title = \juliasymbol, breakable]
    \begin{juliaprompt}
      julia> Z = haagerup_H3_center();
    \end{juliaprompt}
\end{tcolorbox}
Note that loading the category takes about 30 minutes due to the enormous amount of data it contains.
The $\mathcal{Z}(\HI)$ category, obtained this way, is a data structure that allows access to its properties via dedicated functions shown in Table \ref{tab:usefulfuncs}. These are documented in the package documentation \cite{Maeurer_TensorCategories_jl}.
\begin{table}[ht]
  \begin{tabular}{ll}
    \hline
    \textbf{Method} & \textbf{Return Value} \\
    \hline
    \texttt{haagerup\_H3()} & Fusion category $\mathcal{H}_3$ \\
    \texttt{haagerup\_H3\_center()} & Drinfeld center $\mathcal{Z}(\mathcal{H}_3)$ \\
    \texttt{multiplication\_table(C::Category)} & Fusion matrices of the fusion ring of a category \\
    \texttt{smatrix(C::Category)} & $S$-matrix of the category \\
    \texttt{tmatrix(C::Category)} & $T$-matrix of the category \\
    \texttt{fpdim(C::Category)} & Frobenius--Perron dimension of the category \\
    \texttt{simples(C::Category)} & Simple objects of the category \\
    \texttt{F\_symbols(C::Category)} & $F$-symbols of the category \\
    \texttt{R\_symbols(C::Category)} & $R$-symbols of the category \\
    \texttt{numeric\_F\_symbols(C::SixJCategory)} & $F$-symbols as complex numbers \\
    \texttt{numeric\_R\_symbols(C::SixJCategory)} & $R$-symbols as complex numbers\\
    \hline
  \end{tabular}
  \caption{Functions for obtaining data of $\mathcal{H}_3$ and  $\mathcal{Z}(\mathcal{H}_3)$ in the \textsc{TensorCategories.jl} package}
  \label{tab:usefulfuncs}
\end{table}

The following subsections contain examples of their usage and the expected results.
\subsection{Fusion ring}
The Frobenius--Perron dimension can be obtained as follows
\begin{tcolorbox}[title = \juliasymbol, breakable]
    \begin{juliaprompt}
    julia> fpd = fpdim(Z)
	  {a2: 1276.27}
	  julia> minpoly(fpd)
	  x^2 - 1287*x + 13689
    \end{juliaprompt}
\end{tcolorbox}
from which it follows that
\begin{equation}
  D_{FP}(\mathcal{Z}(\mathcal{H}_3)) =\frac{117}{2} \left(3 \sqrt{13}+11\right).
\end{equation}
To obtain the structure constants and assign them to a variable \texttt{mt}, the following code can be used.
\begin{tcolorbox}[title = \juliasymbol, breakable]
    \begin{juliaprompt}
      julia> mt = multiplication_table(Z);
    \end{juliaprompt}
\end{tcolorbox}
The variable \texttt{mt} is a $12\times12\times12$ array whose values are the structure constants of the fusion ring. The code \texttt{mt[}$n$\texttt{,:,:]} returns the $n$'th fusion matrix, i.e., the $n$'th matrix of the left regular representation of the fusion ring. For example, to get the $4$'th fusion matrix, one can evaluate the following code.
\begin{tcolorbox}[title = \juliasymbol, breakable]
    \begin{juliaprompt}
      julia> mt[4,:,:]
      12
       0  0  0  1  0  0  0  0  0  0  0  0
       0  1  1  2  1  1  1  1  1  1  1  1
       0  1  1  1  2  2  1  1  1  1  1  1
       1  2  1  2  1  1  1  1  1  1  1  1
       0  1  2  1  1  2  1  1  1  1  1  1
       0  1  2  1  2  1  1  1  1  1  1  1
       0  1  1  1  1  1  1  1  1  1  1  1
       0  1  1  1  1  1  1  1  1  1  1  1
       0  1  1  1  1  1  1  1  1  1  1  1
       0  1  1  1  1  1  1  1  1  1  1  1
       0  1  1  1  1  1  1  1  1  1  1  1
       0  1  1  1  1  1  1  1  1  1  1  1
    \end{juliaprompt}
\end{tcolorbox}

\subsection{F-symbols, R-symbols, and P-symbols}
The $\mathcal{Z}(\mathcal{H}_3)$ category has a total of 1,630,329 $F$-symbols and 1,401 $R$-symbols. These symbols are stored in dictionaries whose keys are tuples of integers (corresponding to labels of the symbols) and whose values are the values of the symbols. The exact values can be obtained as follows
\begin{tcolorbox}[title = \juliasymbol, breakable]
    \begin{juliaprompt}
      julia> fs = F_symbols(Z)
      Dict{Vector{Int64}, AbsSimpleNumFieldElem} with 1630329 entries:
        [4, 12, 9, 5, 7, 1, 1, 6, 1, 1]   => 611786//765375*_a^47 - 3134042//331662
        [3, 9, 3, 5, 7, 1, 1, 7, 1, 1]    => -962781102528//1842313687825*_a^47 + 1
        [9, 2, 12, 4, 10, 1, 1, 2, 1, 1]  => 14209426039//60020962625*_a^47 - 95649
      julia> rs = R_symbols(Z)
      Dict{Vector{Int64}, AbsSimpleNumFieldElem} with 1401 entries:
        [7, 3, 4, 1, 1]    => -10068//86125*_a^47 - 1234//17225*_a^46 + 341//86125*_
        [11, 12, 3, 1, 1]  => -197//1625*_a^47 - 10457//1625*_a^34 - 231139//1625*_a
        [9, 5, 5, 1, 1]    => -98//1625*_a^41 - 5198//1625*_a^28 - 114741//1625*_a^1
      julia> ps = P_symbols(Z)
      Dict{Vector{Int64}, Nemo.AbsSimpleNumFieldElem} with 12 entries:
        [3]  => 1
        [7]  => 1
        [4]  => 1
        [2]  => 1
        [10] => 1
        [5]  => 1
        [12] => 1
        [8]  => 1
        [1]  => 1
        [6]  => 1
        [11] => 1
        [9]  => 1
    \end{juliaprompt}
\end{tcolorbox}
Note that their entries are not sorted lexicographically on the keys. These exact values are expressed over a number field $\mathbb{Q}(x)$ that can be obtained as follows
\begin{tcolorbox}[title = \juliasymbol, breakable]
    \begin{juliaprompt}
      julia> base_ring(Z)
      Number field with defining polynomial x^48 - x^47 + 2*x^46 - 2*x^45 + 2*x^44 - x^43 - x^42 + 4*x^41 - 8*x^40 + 12*x^39 - 15*x^38 + 15*x^37 - 10*x^36 + 51*x^35 - 31*x^34 + 57*x^33 - 27*x^32 + 2*x^31 + 59*x^30 - 141*x^29 + 229*x^28 - 313*x^27 + 342*x^26 - 285*x^25 + 85*x^24 + 285*x^23 + 342*x^22 + 313*x^21 + 229*x^20 + 141*x^19 + 59*x^18 - 2*x^17 - 27*x^16 - 57*x^15 - 31*x^14 - 51*x^13 - 10*x^12 - 15*x^11 - 15*x^10 - 12*x^9 - 8*x^8 - 4*x^7 - x^6 + x^5 + 2*x^4 + 2*x^3 + 2*x^2 + x + 1
  over rational field
    \end{juliaprompt}
\end{tcolorbox}
The field $\mathbb Q(x)$ obeys $\mathbb{Q}(x)\cong M(\mathcal H_3)(\xi_{39})$ where we use the convention that $\xi_{n}:= \exp\left(\frac{2 \pi i}{n}\right)$, i.e., $\xi_{39}$ is a primitive $39$th root of unity. In particular, the field is a cyclotomic extension of the minimal field of the stored fusion system of $\mathcal{H}_3$. While there might be subfields of $\mathbb{Q}(x)$ over which the $F$-symbols and $R$-symbols are definable, working over cyclotomic extensions has benefits when it comes to computation speed. \\

For applications in physics, it is often desirable to have floating-point approximations of the $F$- and $R$-symbols. These can be obtained as follows:
\begin{tcolorbox}[title = \juliasymbol, breakable]
    \begin{juliaprompt}
      julia> nfs = numeric_F_symbols(Z, precision = 8) # precision is in bits
      Dict{Vector{Int64}, AcbFieldElem} with 1630329 entries:
      [4, 12, 9, 5, 7, 1, 1, 6, 1, 1]   => [-0.049 +/- 6.64e-4] + [0.028 +/- 3.23e-4]*im
      [3, 9, 3, 5, 7, 1, 1, 7, 1, 1]    => [-0.1775 +/- 5.56e-5] + [-0.4611 +/- 4.05e-5]*im
      [9, 2, 12, 4, 10, 1, 1, 2, 1, 1]  => [-0.11520 +/- 7.59e-6] + [-0.7539 +/- 3.68e-5]*im

      julia> nrs = numeric_R_symbols(Z, precision = 8)
      Dict{Vector{Int64}, AcbFieldElem} with 1401 entries:
      [7, 3, 4, 1, 1]    => [0.7091 +/- 6.75e-5] + [0.7050 +/- 7.82e-5]*im
      [11, 12, 3, 1, 1]  => [-0.35 +/- 7.07e-3] + [-0.93 +/- 7.17e-3]*im
      [9, 5, 5, 1, 1]    => [-0.885 +/- 6.32e-4] + [0.465 +/- 4.39e-4]*im
    \end{juliaprompt}
\end{tcolorbox}
The values in the dictionaries are arbitrary-precision complex balls.
Such values can be converted to native floating-point numbers by evaluating, for example,
\begin{tcolorbox}[title = \juliasymbol, breakable]
    \begin{juliaprompt}
      julia> float_rs = Dict(key => convert(ComplexF64, val) for (key, val) in nrs)
      Dict{Vector{Int64}, ComplexF64} with 1401 entries:
      [7, 3, 4, 1, 1]    => 0.709142+0.705066im
      [11, 12, 3, 1, 1]  => -0.354605-0.935016im
    \end{juliaprompt}
\end{tcolorbox}

\begin{remark}
  The $F$-symbols and $R$-symbols are not in a unitary gauge. This is because it was necessary to choose a small field for the $F$-symbols to keep computations reasonable. The unitary $F$-symbols would require a field of degree $16$ \cite{osborne2019f,huang2020f} or $8$ \cite{anyonwiki} and therefore the center would need to be defined over a field of degree at least $192$ or $96$ respectively. This would make the computations infeasible at the moment.
\end{remark}

\subsection{Modular data}
The $S$-matrix of $\mathcal{Z}(\mathcal{H}_3)$ was computed first in \cite{Izumi:2001mi} and simplified in \cite{evans2006modular}. We are able to simplify the $S$-matrix further to

\setcounter{MaxMatrixCols}{20}

\begin{align}\scriptsize
    \begin{pmatrix}
1 &  d^2 &  1 + d^2 &  1 + d^2 &  1 + d^2 &  1 + d^2 &  3d &  3d &  3d &  3d &  3d &  3d \\

 d^2 &  1 &  1 + d^2 &  1 + d^2 &  1 + d^2 &  1 + d^2 &  -3d &  -3d &  -3d &  -3d &  -3d &  -3d \\

 1 + d^2 &  1 + d^2 &  2 + 2d^2 &  -1-d^2 &  -1-d^2 &  -1-d^2 &  0 &  0 &  0 &  0 &  0 &  0 \\

 1 + d^2 &  1 + d^2 &  -1-d^2 &  2 + 2d^2 &  -1-d^2 &  -1-d^2 &  0 &  0 &  0 &  0 &  0 &  0 \\

 1 + d^2 &  1 + d^2 &  -1-d^2 &  -1-d^2 &  -1-d^2 &  2 + 2d^2 &  0 &  0 &  0 &  0 &  0 &  0 \\

 1 + d^2 &  1 + d^2 &  -1-d^2 &  -1-d^2 &  2 + 2d^2 &  -1-d^2 &  0 &  0 &  0 &  0 &  0 &  0 \\

 3d &  -3d &  0 &  0 &  0 &  0 &  a_1 &  a_2 &  a_3 &  a_4 &  a_5 &  a_6 \\

 3d &  -3d &  0 &  0 &  0 &  0 &  a_2 &  a_1 &  a_4 &  a_3 &  a_6 &  a_5 \\

 3d &  -3d &  0 &  0 &  0 &  0 &  a_3 &  a_4 &  a_5 &  a_6 &  a_2 &  a_1 \\

 3d &  -3d &  0 &  0 &  0 &  0 &  a_4 &  a_3 &  a_6 &  a_5 &  a_1 &  a_2 \\

 3d &  -3d &  0 &  0 &  0 &  0 &  a_5 &  a_6 &  a_2 &  a_1 &  a_4 &  a_3 \\

 3d &  -3d &  0 &  0 &  0 &  0 &  a_6 &  a_5 &  a_1 &  a_2 &  a_3 &  a_4
    \end{pmatrix}.
\end{align}

where

\begin{align}
  \begin{aligned}
  d &= \frac{3 + \sqrt{13}}{2} = -\xi_{39}^{23} - \xi_{39}^{17} - \xi_{39}^{14} + \xi_{39}^{12} - \xi_{39}^{10} + \xi_{39}^{9} - \xi_{39}^{4} + \xi_{39}^{3} - \xi_{39} + 2\\
  a_1 &= 3\xi_{39}^{23} - 6\xi_{39}^{21} + 3\xi_{39}^{20} - 6\xi_{39}^{18} + 6\xi_{39}^{17} + 3\xi_{39}^{14} - 3\xi_{39}^{12} + 3\xi_{39}^{10} - 6\xi_{39}^{9} + 3\xi_{39}^{7}  - 3\xi_{39}^{6} \\ & \hspace{1cm} + 6\xi_{39}^{4} - 3\xi_{39}^{3} + 3\xi_{39} - 6 \\
a_2 &= 3\xi_{39}^{23} + 6\xi_{39}^{20} - 3\xi_{39}^{14} + 3\xi_{39}^{12} + 3\xi_{39}^{10} + 6\xi_{39}^{7} - 6\xi_{39}^{6} - 3\xi_{39}^{3} - 3\xi_{39} + 3 \\
a_3 &= -3\xi_{39}^{23} - 6\xi_{39}^{21} - 6\xi_{39}^{18} + 3\xi_{39}^{17} - 3\xi_{39}^{10} - 3\xi_{39}^{9} + 3\xi_{39}^{4} + 3\xi_{39}^{3} + 3 \\
a_4 &= -3\xi_{39}^{23} + 3\xi_{39}^{21} - 6\xi_{39}^{20} + 3\xi_{39}^{18} - 3\xi_{39}^{17} - 3\xi_{39}^{10} + 3\xi_{39}^{9} - 6\xi_{39}^{7} + 6\xi_{39}^{6} - 3\xi_{39}^{4} + 3\xi_{39}^{3} \\
a_5 &= 3\xi_{39}^{23} + 3\xi_{39}^{21} + 3\xi_{39}^{20} + 3\xi_{39}^{18} + 3\xi_{39}^{10} + 3\xi_{39}^{7} - 3\xi_{39}^{6} - 3\xi_{39}^{3} - 3 \\
a_6 &= -6\xi_{39}^{23} + 6\xi_{39}^{21} - 6\xi_{39}^{20} + 6\xi_{39}^{18} - 9\xi_{39}^{17} - 3\xi_{39}^{14} + 3\xi_{39}^{12} - 6\xi_{39}^{10} + 9\xi_{39}^{9} - 6\xi_{39}^{7} + 6\xi_{39}^{6} \\ &\hspace{1cm} - 9\xi_{39}^{4} + 6\xi_{39}^{3} - 3\xi_{39} + 9
  \end{aligned}
\end{align}

And the diagonal of the $T$-matrix is given by

\begin{align}
  (1,1,1,1,\xi_3^{-1}, \xi_3, \xi_{13} ^{-2}, \xi_{13}^{-6}, \xi_{13}^{2}, \xi_{13}^{-5}, \xi_{13}^{6}, \xi_{13}^{5}).
\end{align}

\begin{tcolorbox}[title = \juliasymbol, breakable]
  \begin{juliaprompt}
    julia> Float64.(real.(numeric_smatrix(Z, precision = 8)))
    12
      1.0     10.9082   11.9082   11.9082   11.9082  
     10.9082   1.0      11.9082   11.9082   11.9082       -9.9082    -9.9082    -9.9082
     11.9082  11.9082   23.8164  -11.9082  -11.9082        0.0        0.0        0.0
     11.9082  11.9082  -11.9082   23.8164  -11.9082        0.0        0.0        0.0
     11.9082  11.9082  -11.9082  -11.9082  -11.9082        0.0        0.0        0.0
     11.9082  11.9082  -11.9082  -11.9082   23.8164  
      9.9082  -9.9082    0.0       0.0       0.0         -17.5459   -11.2568    14.832
      9.9082  -9.9082    0.0       0.0       0.0           7.02686   14.832     19.2402
      9.9082  -9.9082    0.0       0.0       0.0          14.832     -2.38843  -17.5459
      9.9082  -9.9082    0.0       0.0       0.0          -2.38843   19.2402   -11.2568
      9.9082  -9.9082    0.0       0.0       0.0     
      9.9082  -9.9082    0.0       0.0       0.0         -11.2568     7.02686   -2.38843

      julia> numeric_twists(Z, precision = 16)
      12-element Vector{AcbFieldElem}:
       1.00000
       1.00000
       1.00000
       1.00000
       [-0.50000 +/- 7.63e-6] + [-0.86602 +/- 7.84e-6]*im
       [-0.500000 +/- 6e-11] + [0.86602 +/- 7.84e-6]*im
       [0.568064746 +/- 9.13e-10] + [-0.822983866 +/- 8.40e-10]*im
       [-0.970941817 +/- 5.06e-10] + [-0.2393156643 +/- 5.06e-11]*im
       [0.568064747 +/- 4.37e-10] + [0.822983866 +/- 2.58e-10]*im
       [-0.748510748 +/- 3.82e-10] + [-0.663122658 +/- 3.86e-10]*im
       [-0.97094 +/- 8.00e-6] + [0.23932 +/- 6.88e-6]*im
       [-0.748510748 +/- 2.65e-10] + [0.663122658 +/- 2.69e-10]*im
  \end{juliaprompt}
\end{tcolorbox} 

\subsection{Correctness of the Data}

The multiplication table, forgetful functor, and $S$-matrix we obtained match those from \cite{Izumi:2001mi, exotic-fusion}. Moreover, the $F$-symbols and $R$-symbols were shown to satisfy the pentagon and hexagon equations by using the \textsc{TensorCategories.jl} package.

\section{The minimal field of a multiplicity-free pivotal fusion category}\label{sec:fields}
As mentioned in \ref{sec:fuscat}, the values of the $F$-, ($R$-)and $P$-symbols of a pivotal (braided) fusion category, and thus also the field to which they belong, depend on a choice of basis (also called a gauge choice) for the multiplicity spaces $H^k_{i,j}$. The performance of algorithms using symbolic algebra depends greatly on the field in which the input is expressed. It is therefore interesting to find a gauge for which the values of the $F$-, $R$-and $P$-symbols lie in a \emph{nice} field.
Ideally, one would like to have a cyclotomic field to which the data belongs. However, it turns out that this is impossible for $\mathcal{H}_2$ and $\mathcal{H}_3$ \cite{morrison2012non}.
The question thus arises to what extent one can \emph{simplify} the field to which the data belongs. It turns out that for multiplicity-free fusion systems over a field extension of $\mathbb{Q}$, there is a notion of a \emph{minimal field}.
Since the $\mathcal{H}_3$ category is not braided, we will only define the notion of a minimal field for a pivotal fusion category, but its extension to a pivotal braided fusion category is straightforward.

Before defining and presenting an algorithm for constructing a minimal field, we will introduce some notation and a few definitions.
We will distinguish between formal symbols and their values. Let $\Phi$ be a pivotal fusion system over $\bbk \supset \mathbb{Q}$. Formal $F$-symbols, $P$-symbols, and gauge-symbols are regarded as functions  $F: \mathbf{L}^6 \to \bbk$, $P:\mathbf{L} \to \bbk$, and $g:\mathbf{L}^3\to \bbk^\times$ whose values are irrelevant. We will denote the lists of formal $F$-symbols, $P$-symbols, and gauge-symbols of a multiplicity-free pivotal fusion system $\Phi$ respectively by $L_F^\Phi$, $L_P^\Phi$, and $L_g^\Phi$, and the ring of rational functions in formal $F$-symbols,
$P$-symbols, and gauge-symbols of $\Phi$ will be denoted by $\mathcal{L}_\Phi = \bbk(L_F^\Phi,L_P^\Phi,L_g^\Phi)$. Note that this ring only depends on $\mathbf{N}$ and not on any of the values of the symbols. We will use the symbol $\gamma$ to denote a formal gauge-transform that transforms formal $F$-and $P$-symbols as in equation \eqref{eq:gt}, where the $g_{i,j}^k$ are instead formal gauge-symbols.

Given a gauge-transform $\varphi$ with values $\Gamma = \gaugespace$, we define the evaluation homomorphisms
\begin{itemize}
  \item $\Phi:\mathcal{L}_\Phi \to \mathcal{L}_\Phi: \rho \mapsto \Phi(\rho)$ as the map which evaluates all occurrences of formal $F$-and $P$-symbols at their values from $\Phi$, and
  \item $\Gamma:\mathcal{L}_\Phi \to \mathcal{L}_\Phi: \rho \mapsto \Gamma(\rho)$ as the map which evaluates all occurrences of formal gauge-symbols at their values from $\Gamma$.
\end{itemize}

The key to finding a gauge with a minimal field for the $F$-and $P$-symbols lies in the construction of a so-called \emph{gauge-split basis} for $\Phi$ \cite{vercleyen2025tables}. In order to define such a basis, we need to distinguish between different ways in which rational expressions of $F$-and $P$-symbols can be gauge-invariant.
\begin{definition} \label{def:inv}
  Let $\Phi$ be a multiplicity-free pivotal fusion system. An element $\rho \in \mathcal{L}_\Phi $ is called
	\begin{itemize}
		\item \emph{de jure gauge-invariant} if it only contains formal $F$-and $P$-symbols and $\gamma(\rho) = \rho$.
		\item \emph{de facto gauge-invariant} if it only contains formal $F$-and $P$-symbols, none of the numerators of $\rho$ evaluate to $0$, and $\Phi(\gamma(\rho)) = \Phi(\rho)$.
	\end{itemize}
\end{definition}
Note that it is not necessary to know the values of the $F$-symbols or $P$-symbols to determine whether $\rho \in \mathcal{L}_\Phi$ is de jure gauge-invariant. Knowledge of $\mathbf{N}$ suffices to figure this out. To determine whether $\rho$ is de facto gauge-invariant, however, extra information about the values of the $F$-and $P$-symbols might be required. Which information depends on the form of $\rho$. If, for example, $\rho$ is a rational monomial in formal $F$-and $P$-symbols, we only need to know which $F$-symbols evaluate to $0$ since the $P$-symbols are always non-zero.

A gauge-split basis of a multiplicity-free pivotal fusion system $\Phi$ can then be defined as follows.
\begin{definition}\label{def:gsb} Let $\mathcal{M}_{\Phi} \subset \mathcal{L}_\Phi$ be the (free abelian) subgroup of monic rational monomials (or words) in formal $F$-, $P$-and gauge-symbols. A \emph{gauge-split-basis} (GSB) of $\Phi$  is a tuple $( I, D )$ where $I = (\iota_1,\ldots,\iota_m) \in \mathcal{M}^m_\Phi$, $D=(\delta_1,\ldots,\delta_n)\in \mathcal{M}^n_\Phi$,
  where $m,n\in \mathbb{N}$ satisfy $m+n = |L_F^\Phi| +|L_P^\Phi|$ and $ n \leq |\mathbf{N}|$,  and
  for which the following properties hold.
	\begin{enumerate}
		\item Each word $\rho$ in $I$ and $D$ only contains formal $F$-and $P$-symbols and $\Phi(\rho)$ is well-defined.
		\item Each word in $I$ is de facto gauge-invariant.
		\item Each word in $D$ is gauge-dependent and moreover for any $n$-tuple $w\in \left(\bbk^\times\right)^n$, there exists a gauge-transform with variables $\Gamma$ such that $\Phi(\Gamma(\gamma(D))) = w$.
		\item For any word $\rho\in\mathcal{M}_\Phi$ there exist $m+n$ unique integers $a_1,\ldots,a_m, b_1,\ldots, b_n$ such that
      \begin{eqnarray}
	      \rho = \prod_{i = 1}^m \iota_i^{a_i} \prod_{i = 1}^n \delta_i^{b_i}\label{eq:uniqueness}
      \end{eqnarray}
	\end{enumerate}
\end{definition}
\begin{remark}
  The term basis comes from property $(4)$, which tells us that any word in formal $F$- and $P$-symbols can be written as a unique multiplication of powers of elements of $I$ and $D$.
\end{remark}
\begin{remark}
  In \cite{hagge2015geometricinvariantsfusioncategories} and the appendix of \cite{davidovichArithmeticModularCategories2013}, a method to check gauge-equivalence of two sets of $F$-symbols by using a gauge-invariant set is discussed. This set is the same as the set $I$ in this paper. While we are certain that the authors of \cite{hagge2015geometricinvariantsfusioncategories} are aware of all the properties of a gauge-split basis, they never defined a term for it.
\end{remark}

\begin{example}
  Let $\Phi = [\mathrm{Fib}]_{1,1,1}$ (also denoted as $\mathrm{FC}^{2,1,0}_{1,1,1,1}$ on the AnyonWiki \cite{anyonwiki}) be the fusion system with $\mathbf{L} = \{1,2\}$ and fusion rules
  \begin{align}
    1 \otimes 1 & \cong 1,\\
    1 \otimes 2 &\cong 2 \otimes 1 \cong 2,\\
    2\otimes 2  &\cong 1 \oplus 2.
  \end{align}
  The following is a gauge split basis for $\Phi$
  \begin{IEEEeqnarray*}{rcl}
    I &= \left( \vphantom{\frac{\Fs{2,1,2,1,2,2}}{\Fs{1,1,2,2,1,2} \Fs{2,1,1,2,2,1}}}\right. &
    \Fs{1,1,1,1,1,1},\Fs{1,2,1,2,2,2},P_1,\Fs{1,1,2,2,1,2} \Fs{1,2,2,1,2,1},\Fs{2,1,1,2,2,1}
       \Fs{2,2,1,1,1,2},\frac{\Fs{2,1,2,1,2,2}}{\Fs{1,1,2,2,1,2} \Fs{2,1,1,2,2,1}},\\
       & &\Fs{1,1,2,2,1,2} \Fs{2,1,1,2,2,1}
       \left. \Fs{2,2,2,2,1,1},P_2 \Fs{1,1,2,2,1,2} \Fs{2,1,1,2,2,1}
       \vphantom{\frac{\Fs{2,1,2,1,2,2}}{\Fs{1,1,2,2,1,2} \Fs{2,1,1,2,2,1}}}\right),\\
    D &=& \left(\Fs{1,1,2,2,1,2},\Fs{2,1,1,2,2,1}\right).
  \end{IEEEeqnarray*}
  Any element of $I$ is gauge-invariant, and each element of $D$ can be set to any non-zero value by choosing a suitable gauge transform. Moreover, any $F$- and $P$-symbol of $\Phi$ can be written in a unique way as a multiplication of powers of elements of $I$ and $D$.
\end{example}

We have the following useful result
\begin{lemma}\label{lem:GSB}
  Any multiplicity-free pivotal fusion system has an infinite number of gauge-split bases.
\end{lemma}
Before proving Lemma \ref{lem:GSB}, we will introduce a useful definition and lemma. Consider the free abelian group $\fpmod$, obtained by taking the direct product of $|\lfpp|$ copies of $\mathcal M_\Phi$.
The list $\lfpp$ of all formal $F$- and $P$-symbols, with length $|\lfpp|$, can then be regarded as an element of $\fpmod$ and set of $m$ words in formal $F$-and $P$-symbols can be generated using the following product
\begin{definition}\label{def:powerdot}
   Let $\Phi = (\mathbf{L},*,\mathbf{N},\mathbf{F},\mathbf{P})$ be a multiplicity-free pivotal fusion system and $m\in \mathbb{N}$ and $n:=|\lfpp|$. We define
   \begin{eqnarray}
     -\uparrow-&:& \mathrm{Mat}_{m\times n}(\mathbb{Z}) \times \fpmod \to \mathcal{M}_\Phi^m,\nonumber\\
     && [M \uparrow v]^i = \prod_{j=1}^{n} \left([v]^j\right)^{[M]^i_j},
   \end{eqnarray}
   where $[v]^j$ is the $j$'th element of $v$ and $[M]^i_j$ is the $(i,j)$'th matrix coefficient of $M$.
\end{definition}
\begin{remark}
  If the multiplication of elements of the abelian group $\fpmod$ were written additively, then we could just use the normal matrix product notation. We do not do this to avoid confusion.
\end{remark}
\begin{lemma}\label{lem:module}
  If $A\in \mathrm{Mat}_{m\times n}(\mathbb{Z}), B\in \mathrm{Mat}_{n\times l}(\mathbb{Z})$ then  $A\uparrow (B \uparrow v) = (AB) \uparrow v$. In particular, $\uparrow$ turns $\fpmod$ into a left $\mathrm{Mat}_{|\lfpp|\times|\lfpp|}(\mathbb{Z})$-module.
\end{lemma}
We will now prove Lemma \ref{lem:GSB}. While it is possible to provide a short proof by using the language of free groups generated by formal $F$- and $P$ symbols, we decided to give a longer, more explicit proof via the construction of the GSB. This has the advantage of being easy to translate to an algorithm such as the one that Anyonica uses.
\begin{proof}
Let $\Phi = (\mathbf{L},*,\mathbf{N},\mathbf{F}, \mathbf{P})$ be a multiplicity-free pivotal fusion system over $\bbk \supset \mathbb{Q}$. Let $n_z$ be the number of $F$-symbols of $\Phi$ whose value is $0$. Let $L$ be $L_{FP}$ without the formal $F$-symbols that evaluate to $0$, let $k=|\lfpp|-n_z$ be its length, and let $l = |L_F^\Phi|-n_z$ be the number of formal $F$-symbols that do not evaluate to $0$. Let $\gamma$ be a formal gauge transform and define
\begin{eqnarray}
  G_i :=
  \left\{
  \begin{array}{lll}
    \frac{\g{a_i,b_i,e_i}\g{e_i,c_i,d_i}}{\g{a_i,f_i,d_i}\g{b_i,c_i,f_i}} & \text{ if } & i \leq l \\
    \frac{\g{a_i^*,a_i,1}\g{1,a_i^*,a_i^*}}{\g{a_i^*,1,a_i^*}\g{a_i,a_i^*,1}} & \text{ if } & i > l
  \end{array}\right.
\end{eqnarray}
as the overall multiplicative gauge factor that appears in the gauge transform $[L]^i \mapsto \gamma([L]^i) = G_{i} [L]^i$
Any word in the formal non-zero $F$-and $P$-symbols $[L]^1, \ldots, [L]^k $ can be written as
  \begin{eqnarray}
    \rho_p &=& p \uparrow L \\
    &=& ([L]^1)^{[p]_1} ([L]^2)^{[p]_2} \cdots ([L]^k)^{[p]_k}  ,\quad [p]_i \in \mathbb{Z},
  \end{eqnarray}
  where we have identified $p \in \mathrm{Mat}_{1\times k}(\mathbb Z)$ with an integer row vector.
For $\rho_{p}$ to be de facto gauge-invariant, we need to have that
\begin{IEEEeqnarray}{rrCl}
  &\Phi(\gamma(\rho_{p})) &=& \Phi(\rho_{p}) \\
  \Leftrightarrow & \Phi((G_1 [L]^1)^{[p]_1} \! \cdots ( G_k [L]^k)^{[p]_k}) &=& \Phi(([L]^1)^{[p]_1} ([L]^2)^{[p]_2} \cdots ([L]^k)^{[p]_k}) \label{eq:1} \\
  \Leftrightarrow  & \left(G_1\right)^{[p]_1}  \cdots \left(G_k\right)^{[p]_k} &=& 1\label{eq:2}, \\
\end{IEEEeqnarray}
where equations \eqref{eq:1} and \eqref{eq:2} are equivalent because none of the elements of $L$ evaluate to $0$.
After expanding the $G_i$ equation \eqref{eq:2} is equivalent to
\begin{IEEEeqnarray}{rCl}
  \left(\!\frac{\g{a_1,b_1,e_1}\g{e_1,c_1,d_1}}{\g{a_1,f_1,d_1}\g{b_1,c_1,f_1}}\!\right)^{[p]_1}\!\!\!\!\!\! \cdots \left(\!\frac{\g{a_l,b_l,e_l}\g{e_l,c_l,d_l}}{\g{a_l,f_l,d_l}\g{b_l,c_l,f_l}}\!\right)^{[p]_l}
  \left(\! \frac{\g{a_{l+1}^*,a_{l+1},1}\g{1,a_{l+1}^*,a_{l+1}^*}}{\g{a_{l+1}^*,1,a_{l+1}^*}\g{a_{l+1},a_{l+1}^*,1}}\!\right)^{[p]_{l+1}}\!\!\!\!\!\! \cdots
  \left(\!\frac{\g{a_{k}^*,a_{k},1}\g{1,a_{k}^*,a_{k}^*}}{\g{a_{k}^*,1,a_{k}^*}\g{a_{k},a_{k}^*,1}}\!\right)^{[p]_k}
  &=& 1. \label{eq:gauge_invariance}
\end{IEEEeqnarray}
Now assume that the $\g{a,b,c}$ are ordered lexicographically on $(a,b,c)$ in a list of length $h$ and write $ g_j$ for the $j$th variable that appears in this list. By collecting factors $g_j$ in equation (\ref{eq:gauge_invariance}) we get that
\begin{eqnarray}\label{eq:gauge_invariance2}
  g_1^{\sum_{i_1}[M]^1_{i_1}  [p]_{i_1}} g_2^{\sum_{i_2} [M]^2_{i_2} [p]_{i_2}} \cdots g_h^{\sum_{i_h} [M]^{l}_{i_h} [p]_{i_h}} = 1, 
\end{eqnarray}
where $M$ is a matrix with integer coefficients. Equation (\ref{eq:gauge_invariance2}) can only be satisfied if each exponent of each $g_i$ is identically zero, i.e.,
\begin{eqnarray}\label{eq:matrix_z}
  M.p^T = \mathbb{0}.
\end{eqnarray}
The space of vectors $p$ for which $\rho_p$ is gauge-invariant is isomorphic to the kernel of $M$, while the span of the vectors $q$ for which $\rho_q$ is gauge-dependent is isomorphic to the co-kernel of $M$. The kernel and co-kernel of $M$ can be constructed via a Smith decomposition. Indeed, if we write $M = U.S.V$, with $U, V$ invertible integer matrices and $S$ diagonal, then $m = \dim\left( \ker M\right)$
equals the number of zero columns of $S$ and $\ker M$ is spanned by the last $m$
columns of $V^{-1}$. The first $k-m =: n$ columns of $V^{-1}$ span $\coker M$.
Note that any integer matrix has a Smith decomposition.

Let $A =  (V^{-1})^T \oplus \mathbb{1}_{n_Z\times n_Z}$. Let $K$ be the concatenation of $L$ with the list of formal $F$-symbols, in any order, that evaluate to $0$. In particular $|K| = |\lfpp|$ and the last $n_z$ elements of $K$ are formal $F$-symbols that evaluate to $0$.
Let $v = A \uparrow K$, $D \in \mathcal{M}_\Phi^n$ be the vector consisting of the first $n$ elements of $v$, and $I \in \mathcal{M}_\Phi^{|\mathbf{F}|-n}$ be the vector of the remaining elements of $v$.
  $(I,D)$ then has the following properties.
\begin{enumerate}
  \item Every word $\rho$ in $I$ and $D$ is a word in formal $F$-and $P$-symbols from $\lfpp$. Moreover, since neither $I$ nor $D$ contains negative powers of $F$-symbols that evaluate to $0$, we have that $\Phi(\rho)$ is well-defined for every word $\rho$ in $I$ and $D$.
  \item By construction, every element of $I$ is de facto gauge-invariant.
  \item By construction, every element of $D$ is gauge dependent. Moreover, let $w \in (\bbk^\times)^n$, then
  \begin{eqnarray}
    \Phi(\Gamma(\gamma(D))) = w
  \end{eqnarray}
  has at least one solution for the variables $\gaugespace$. Indeed, let $B = (V^{-1})^T$, $G$ be the vector with coefficients $[G]^i = G_i$, $g$ be the vector with coefficients $[g]^i=g_i$, then for all $j=1,\ldots,n$

  \begin{eqnarray}
    [\Phi(\Gamma(\gamma(D)))]^j &=& \Phi(\Gamma(\gamma ( [ A \uparrow K ]^j ) ) ) =  \Phi(\Gamma(\gamma ( [ B \uparrow L ]^j ) ) )\\
    &=& \prod_{i=1}^{|L|}\Phi(\Gamma( \gamma([L]^i ) ) )^{[B]^j_i}
    = \prod_{i=1}^{|L|}\Gamma( [G]^i )^{[B]^j_i} \prod_{i=1}^{|L|}\Phi([L]^i)^{[B]^j_i}\\
    &=& [ B \uparrow \Gamma(G) ]^j [ B \uparrow \Phi(L) ]^j.
\end{eqnarray}
Since
\begin{eqnarray}
  [ B \uparrow \Gamma(G) ]^j = [ B \uparrow M^T \uparrow \Gamma(g) ]^j =  [ (B M^T)  \uparrow \Gamma(g) ]^j = [ (US)^T  \uparrow \Gamma(g) ]^j,
\end{eqnarray}
we have that
\begin{eqnarray}
    \Phi(\Gamma(\gamma(D))) &=& w \\
      \Leftrightarrow [ (US)^T  \uparrow \Gamma(g) ]^j &=&  \frac{[w]^j}{[ B \uparrow \Phi(L) ]^j} =: [X]^j,\quad j = 1, \ldots, n
\end{eqnarray}
  where we used the fact that none of the elements of $L$ evaluate to $0$.
  Let $C$ be the $n \times h$ matrix with entries $[C]^i_j = [US]^j_i$ then $[ (US)^T  \uparrow \Gamma(g) ]^j = [ C \uparrow \Gamma(g)]^j$ for $j=1,\ldots,n$. Since $C$ has rank $n$, there exists a $\tilde{C}$ such that $\tilde{C}C = \mathbb{1}_{n\times n}$ and we find that for such a $\tilde{C}$
  \begin{eqnarray}
    [\Gamma(g)]^i  &=&   [ \tilde{C} \uparrow X]^i, \quad i = 1,\ldots h\\
    \Rightarrow \Phi(\Gamma(\gamma(D))) &=& w
  \end{eqnarray}
  \item Any word $\rho$ is of the form $\rho = p \uparrow L_F$ for a unique $p$. Since $A$ is invertible $D\oplus I = v = A \uparrow L_F$. Because of Lemma \ref{lem:module} $ L_F = A^{-1}\uparrow (D\oplus I)$ and thus $\rho = (p.A^{-1})\uparrow (D\oplus I)$ which is of the form presented in equation (\ref{eq:uniqueness}).
\end{enumerate}
In \cite{jager2005reduction}, it is shown that there are an infinite number of integer transformation matrices $U,V$ that result in the same diagonal matrix $S$ and therefore there are an infinite number of different gauge-split bases.
\end{proof}

\begin{theorem}\label{thm:min_field}
Let $\Phi$ be a multiplicity-free pivotal fusion system over $\bbk \supset \mathbb{Q}$ and $(I,D)$ be a GSB for $\Phi$. Let $\lfpp$ be the concatenation of the lists $L_{F}^\Phi$ and $L_{P}^\Phi$ and $M(I)=\mathbb{Q}(\Phi(I))$ be the minimal field extension of $\mathbb{Q}$ that contains all values of $\Phi(I)$. Then
  \begin{enumerate}
    \item there exists a gauge transform such that $\Phi(\Gamma(\gamma(s))) \in M(I)$ for all $s\in \lfpp$, i.e., there exists a gauge such that all $F$-and $P$-symbols belong to the field $M(I)$.
    \item $M(I)$ is independent of the choice of GSB, i.e., $M(I_1) = M(I_2)$ for any two GSBs $(I_1,D_1)$ and $(I_2,D_2)$ of $\Phi$.
    \item For any gauge transform, we have that $M(I)\subseteq \mathbb{Q}(\Phi(\Gamma(\gamma(\lfpp))))$. In other words, for any gauge, the field to which the values of the $F$-and $P$-symbols belong must contain $M(I)$ as a subfield.
    \item $M(I)$ is invariant under the action of fusion ring automorphisms, i.e., for any fusion system $\Phi'$ whose symbols are obtained by applying a fusion ring automorphism to $\Phi$ we have that $\mathbb{Q}(\Phi'(I)) = \mathbb{Q}(\Phi(I)) = M(I)$. In particular, $M(I)$ is an invariant of the fusion category represented by $\Phi$.
  \end{enumerate}
\end{theorem}

\begin{definition}
  The field $M(I)$ in Theorem \ref{thm:min_field} is called the \emph{minimal field} of $\Phi$ (and of its fusion category $\calc$). We will denote it by $M(\Phi)$ or $M(\calc)$. If $M(\Phi)$ is contained in a cyclotomic field, then the cyclotomic field with the lowest degree containing $M(\Phi)$  is called the \emph{minimal cyclotomic field} of $\Phi$.
  $M(\calc)$
\end{definition}

\begin{remark}
Note that the matrix $A$ in the proof of Lemma \ref{lem:GSB} can be used to perform a gauge transform of a set of $F$--and $P$--symbols to a gauge that requires only the minimal field to express these symbols. Indeed, given a gauge-split basis $(I,D)$ of a pivotal fusion system $\Phi$. We can set every element of $D$ equal to $1$ and then obtain the new values of the $F$--and $P$--symbols as follows. Let $n_D$ be the length of $D$ and $n_I$ be the length of $I$, then $L'_{FP} = \Phi\left(A^{-1}\uparrow v\right)$ where $[v]^i = 1$ for $i = 1, \ldots, n_D$
and $[v]^i=[I]^i$ for $i = n_D+1, \ldots, n_D+n_I $, provides a new set of $F$--and $P$--symbols which can be expressed in the field $M(\Phi)$.
\end{remark}

\begin{remark}
  It is possible to define a GSB for a multiplicity-free fusion system, and likewise, it is possible to extend the notion of a GSB to a multiplicity-free (pivotal) braided fusion system. Lemma \ref{lem:GSB} and Theorem \ref{thm:min_field} can be adapted to work for these cases, and the notion of a minimal field stays intact.
\end{remark}

For the categories $\mathcal{H}_2$ and $\mathcal{H}_3$, the fusion systems from the AnyonWiki and Anyonica provide the following minimal fields
\begin{itemize}
  \item $M(\mathcal{H}_2) = \mathbb{Q}(\alpha)$ with $\alpha^4 - \alpha^3 - \alpha^2 - \alpha + 1 = 0$,
  \item $M(\mathcal{H}_3) = \mathbb{Q}(\beta)$  with $\beta^4 + \beta^3 - \beta^2 + \beta + 1 = 0$.
\end{itemize}

\end{document}

%% file: Tikzpictures.tex
\usetikzlibrary{shapes.geometric}
\tikzset{
    baseline = (current bounding box.center)
}

\newcommand{\tzCanonicalBasisLeft}[7]{
    \begin{tikzpicture}[baseline = (current bounding box.center)]
        \draw (0,0) -- (0.5,1) -- node[pos = 0.5, below right = -5pt, scale = 0.7] {$#5$} (1,2) -- (1,3);
        \draw (1,0) -- (0.5,1);
        \draw (2,0) -- (1,2);
        \node[below] at (0,0) {$#1$};
        \node[below] at (1,0) {$#2$};
        \node[below] at (2,0) {$#3$};
        \node[above] at (1,3) {$#4$};
        \node[left] at (0.5,1) {$#6$};
        \node[left] at (1,2) {$#7$};
    \end{tikzpicture}
}

\newcommand{\tzCanonicalBasisRight}[7]{
    \begin{tikzpicture}[baseline = (current bounding box.center)]
        \draw (2,0) -- (1.5,1) -- node[pos = 0.5, below left = -5pt, scale = 0.7] {$#5$} (1,2) -- (1,3);
        \draw (1,0) -- (1.5,1);
        \draw (0,0) -- (1,2);
        \node[below] at (0,0) {$#1$};
        \node[below] at (1,0) {$#2$};
        \node[below] at (2,0) {$#3$};
        \node[above] at (1,3) {$#4$};
        \node[right] at (1.5,1) {$#6$};
        \node[right] at (1,2) {$#7$};
    \end{tikzpicture}
}

\newcommand{\tzBasisVector}[4 ]{
    \begin{tikzpicture}
        \draw (0,0) -- (0.5,1) -- (0.5,2);
        \draw (1,0) -- (0.5,1); 

        \node[below] at (0,0) {$#1$};
        \node[below] at (1,0) {$#2$};
        \node[above] at (0.5,2) {$#3$};
        \node[left] at (0.5,1) {$#4$};
    \end{tikzpicture}
}